\documentclass{article}

\usepackage{amsmath}
\usepackage{amssymb}
\usepackage{latexsym}
\usepackage{amsxtra}
\usepackage{amscd}
\usepackage{theorem}
\usepackage{xypic}

\theoremstyle{plain}

{\theorembodyfont{\slshape}

        \newtheorem{thm}{Theorem}[section]
        \newtheorem{cor}[thm]{Corollary}
        \newtheorem{lem}[thm]{Lemma}
        \newtheorem{prop}[thm]{Proposition}

}

{\theorembodyfont{\rmfamily}

        \newtheorem{defn}[thm]{Definition}
        
        \newtheorem{rem}[thm]{Remark}

        \newtheorem{ass}[thm]{Assumption}

}

\renewcommand{\em}{\sl}

\newcommand{\proof}{{\bf Proof:\ }}
\newcommand{\Endproof}{\hspace*{\fill} $\Box$ \vspace{1ex} \noindent }

\makeatletter
\renewcommand{\subsection}{\@startsection{subsection}{2}%
        {\z@}{-3.25ex plus -1ex minus-.2ex}{-1em}{\bf}}
\makeatother

\newcommand{\NN}{\mathbb{N}}
\newcommand{\ZZ}{\mathbb{Z}}
\newcommand{\QQ}{\mathbb{Q}}

\newcommand{\RR}{\mathbb{R}}

\newcommand{\C}{\mathcal{C}}
\newcommand{\OO}{\mathcal{O}}

\newcommand{\Gal}{{\rm Gal}}

\newcommand{\Br}{{\rm Br}}

\newcommand{\p}{\mathfrak{p}}

\newcommand{\Res}{{\rm Res}}

\newcommand{\inj}{\hookrightarrow}
\newcommand{\To}{\;\longrightarrow\;}
\newcommand{\iso}{\stackrel{\sim}{\to}}

\newcommand{\lpfeil}[1]{\stackrel{#1}{\To}}

\newcommand{\abs}[1]{\lvert#1\rvert}

\newcommand{\sw}{{\rm sw}}
\newcommand{\dsw}{{\rm dsw}}
\newcommand{\rsw}{{\rm rsw}}

\newcommand{\kb}{\bar{k}}
\newcommand{\Kb}{\bar{K}}
\newcommand{\Mb}{\bar{M}}
\newcommand{\Lb}{\bar{L}}

\newcommand{\xb}{\bar{x}}
\newcommand{\yb}{\bar{y}}
\newcommand{\zb}{\bar{z}}
\newcommand{\ub}{\bar{u}}
\newcommand{\vb}{\bar{v}}
\newcommand{\wb}{\bar{w}}
\newcommand{\ab}{\bar{a}}
\newcommand{\bb}{\bar{b}}
\newcommand{\cb}{\bar{c}}
\newcommand{\db}{\bar{d}}
\newcommand{\chib}{\bar{\chi}}
\newcommand{\deltab}{\bar{\delta}}
\newcommand{\omegab}{\bar{\omega}}

\newcommand{\deltat}{\tilde{\delta}}

\begin{document}

\title{Fiercely ramified cyclic extensions of $p$-adic fields with imperfect
  residue field}

\author{Stefan Wewers\\ Institut f\"ur Reine Mathematik\\
          Universit\"at Ulm}

\date{}

\maketitle

\begin{abstract}
  We study the ramification of fierce cyclic Galois extensions of a \mbox{local}
  field $K$ of characteristic zero with a one-dimensional residue field of
  characteristic $p>0$. Using Kato's theory of the refined Swan conductor, we
  associate to such an extension a {\em ramification datum}, consisting of a
  sequence of pairs $(\delta_i,\omega_i)$, where $\delta_i$ is a positive
  rational number and $\omega_i$ a differential form on the residue field of
  $K$. Our main result gives necessary and sufficient conditions on such
  sequences to occur as a ramification datum of a fierce cyclic extension of
  $K$.

  \vspace{1ex}\noindent
\end{abstract}

\section{Introduction}

\subsection{A classical result on wildly ramified cyclic extensions}

Let $K$ be a complete discrete valued field with residue field $\Kb$ of
characteristic $p>0$. Let $L/K$ be a finite Galois extension with Galois group
$G$. Let us assume, for the moment, that the residue field $\Kb$ is
perfect. Then the classical theory of ramification groups gives rise to a
descending filtration $(G^t)_{t\in\QQ, t\geq -1}$ on $G$ (the {\em upper numbering
  filtration}), which has the following nice properties (see \cite{SerreCL}):
\begin{itemize}
\item
  The definition of the filtration is compatible with taking quotients of $G$.
\item The subgroup $G^0$ is the inertia subgroup of $G$. Moreover, the graded
  pieces ${\rm Gr}^t(G)$ for $t>0$ are elementary abelian and killed by $p$.
\item If $G$ is abelian then all the breaks of the filtration are integers
  (Theorem of Hasse and Arf).
\end{itemize}

Now suppose in addition that the field $K$ has characteristic $p$,
$G\cong\ZZ/p^n\ZZ$ is cyclic of $p$-power order and that $L/K$ is totally
ramified. It is clear from the above that the sequence of breaks of the
filtration $(G^t)_t$ is now of the form
\[
   0<u_1<u_2<\ldots<u_n,
\]
where each $u_i$ is a positive integer such that $\abs{G^{u_i}}=p^{n-i+1}$. In
this situation we have the following (more or less well known) result.

\begin{thm} \label{introthm1}
  Suppose that ${\rm char}(K)=p$. 
  Let $0<u_1<u_2\ldots<u_n$ be a strictly increasing sequence of
  integers. Then $(u_i)$ is the sequence of ramification breaks associated to
  a totally ramified Galois extension $L/K$ with cyclic Galois group of order
  $p^n$ if and only if the following holds: 
  \begin{enumerate}
  \item
    $p\nmid u_1$, 
  \item
    for $i>1$ we either have $u_i=pu_{i-1}$ or $u_i>pu_{i-1}$ and $p\nmid
    u_i$. 
  \end{enumerate}
\end{thm}

A more general version of this result is proved in \cite{ObusPries10}.  As
stated above, the theorem is an immediate consequence of the classical formula
of H.L.\ Schmid (\cite{HLSchmid37}, Satz 3) which computes the breaks $u_i$ in
terms of the Artin-Schreier-Witt representation of $L/K$. See also
\cite{Thomas05} and \cite{Garuti99}.
   
The goal of the present paper is to prove results that are analogous to
Theorem \ref{introthm1} in the case when $K$ has characteristic zero and the
residue field $\Kb$ is `one-dimensional' in a suitable sense (in particular:
not perfect!). These results are used in joint work with Andrew Obus
(\cite{ObusWewers}) on the problem of `lifting' a cyclic Galois
extension as in Theorem \ref{introthm1} to characteristic zero (this problem
is called the {\em local lifting problem}). In the context of these
applications, the results of the present paper can be seen as a `deformation' of
Theorem \ref{introthm1}. 

\subsection{Fierce extensions}

If the residue field of $K$ is not perfect, the theory of ramification groups
becomes more delicate. A general construction of a filtration
$(G^t)_t$ for the Galois group of an extension of the field $K$ has been given
by Abbes and Saito in \cite{AbbesSaito02}. In subsequent papers, the same
authors have also shown that their filtration has all the expected properties
(e.g.\ the Hasse-Arf-Theorem holds).

From the point of view of ramification theory, the situation considered in the
present paper is rather special. We therefore get by with using the much
simpler construction of ramification groups explored in \cite{KatoVC}, which
is similar to the classical construction. We do use in an essential way the
beautiful results of Kato which gives rise to the {\em refined Swan conductor}
and characterizes it in terms of higher class field theory (\cite{KatoLCFT1},
\cite{KatoDSC}, \cite{KatoVC}).

Our assumptions on $K$ are as follows. We assume that $K$ has characteristic
zero and is complete with respect to a discrete valuation $v$ with residue
field $\Kb$. We normalize $v$ such that $v(p)=1$. The crucial assumption is
that the residue field $\Kb$ has a $p$-basis of length one (i.e.\
$[\Kb:\Kb^p]=p$) and that $H^1(\Kb,\ZZ/p\ZZ)\neq 0$. For instance, this
assumption holds if $\Kb$ is itself a local field with perfect residue field. 
 
A finite separable extension $L/K$ is called {\em fierce} if the extension of
residue fields $\Lb/\Kb$ is purely inseparable and
$[L:K]=[\Lb:\Kb]$.\footnote{Some authors call such an extension {\em
    ferociously ramified}.} Note that a fierce extension is {\em weakly
  unramified}, i.e.\ that $K$ contains a prime element of $L$. Moreover, our
assumption on $\Kb$ implies that any finite fierce extension $L/K$ is
monogenic (which means that the corresponding extension of valuation rings is
generated by one element). So if $L/K$ is fierce and Galois with Galois group
$G$, then one can define a filtration of higher ramification groups $(G^t)_t$
on $G$ essentially as in the classical situation (\cite{KatoVC}).
The point is that for the kind of applications we have in mind (related to
reduction and lifting of covers of curves) it is very natural to replace, if
necessary, the field $K$ by a suitable {\em constant extension} (see Section
\ref{prelim2}, Definition \ref{constantdef}). Using a result of Epp
(\cite{Epp73}) we may therefore assume from the start that our Galois
extension $L/K$ is fierce.

\subsection{Ramification data for fierce cyclic extensions}

Let us now assume that $L/K$ is a fierce Galois extension with cyclic Galois
group $G\cong\ZZ/p^n\ZZ$. Using the results of Kato mentioned before, we
associate to $L/K$ a so called {\em ramification datum}
$(\delta_i,\omega_i)_{i=1,\ldots,n}$. Here $\delta_1,\ldots,\delta_n$ are
positive rational numbers and $\omega_1,\ldots,\omega_n$ are differential
forms on the residue field $\Kb$. For fixed $i$, the pair
$(\delta_i,\omega_i)$ represents the {\em refined Swan conductor} of a
character $\chi_i$ of $G$ of order $p^i$, in the sense of \cite{KatoDSC}. The
rational numbers $\delta_1<\delta_2<\ldots<\delta_n$ are simply the breaks of
the filtration $(G^t)_t$; the differentials $\omega_1,\ldots,\omega_n$ contain
finer information on the ramification of $L/K$. For instance, if
$v_1:\Kb^\times\to\ZZ$ is a normalized discrete valuation on the residue field
of $K$ and $\tilde{v}:\Kb\to\QQ\times\ZZ$ is the rank-$2$-valuation obtained
from the pair $(v,v_1)$, then the breaks of the ramification filtration on $G$
with respect to the valuation $\tilde{v}$ are 
\[
       (\delta_i,v_1(\omega_i)+1), \quad i=1,\ldots,n.
\]
This follows from \cite{KatoVC}, Corollary 4.6.
 
Our main result (Theorem \ref{thm1} and Theorem \ref{thm2}) gives a necessary
and sufficient condition on a tupel $(\delta_i,\omega_i)_{i=1,\ldots,n}$ to
occur as the ramification datum of a fierce cyclic extension $L/K$. The full
statement of this condition is a bit lengthy (see Theorem \ref{thm1}). The
following theorem states a partial result which is analogous to Condition (ii)
in Theorem \ref{introthm1} (the letter $\C$ stands for the Cartier operator).

\begin{thm} \label{introthm2}
  Let $L/K$ be a fierce Galois extension with cyclic Galois group
  $G\cong\ZZ/p^n\ZZ$. Let $(\delta_i,\omega_i)_{i=1,\ldots,n}$ be the
  ramification datum associated to $L/K$. Suppose $\delta_{i-1}<1/(p-1)$ for
  some $i>1$. Then either
  \[
      \delta_i=p\,\delta_{i-1}, \quad \C(\omega_i)=\omega_{i-1},
  \]
  or
  \[
       \delta_i>p\,\delta_{i-1}.
  \]
  In the latter case, we have $\C(\omega_i)=\omega_i$ if and only if
  $\delta_i=p/(p-1)$ and $\C(\omega_i)=0$ otherwise.
\end{thm}

The inequality $\delta_i\geq p\,\delta_{i-1}$ resulting from this theorem (which
holds if $\delta_{i-1}\leq 1/(p-1)$) is not an entirely new result. For
instance, it follows from the results of Hyodo (\cite{Hyodo87}). However,
Hyodo's results are valid in much greater generality then ours and are
therefore weaker. In particular the statements in Theorem \ref{introthm2}
concerning the differential forms $\omega_i$ are false if we omit the
assumption $[\Kb:\Kb^p]=p$. It is precisely these relations beween the
differentials $\omega_i$ which are crucial for the applications of our results
in \cite{ObusWewers}.

\vspace{4ex}\noindent
{\bf Acknowledgments:} I am very grateful for the referee for a careful
reading of a previous version of this paper and for pointing out several
inaccuracies.

\section{Preliminaries} \label{prelim}

We introduce the basic notation and assumptions concerning the field $K$ we
will be working with. We also discuss (almost) constant and fierce extensions
of $K$, following \cite{ZhukovImperfect}.

\subsection{} \label{notation}

Let $K$ be a complete discrete valuation field of characteristic zero, with
residue field $\Kb$ of characteristic $p>0$. The valuation on $K$ is denoted
by $v$, and it is normalized such that $v(p)=1$. Let $e\in\NN$ be the absolute
ramification index of $K$, so that $v(K^\times)=\ZZ e^{-1}$. For $t\in
\ZZ e^{-1}$ we set
\[
     \p_K^t := \{\, x\in K \mid v(x)\geq t \,\}
\]
and
\[
    U_K^t := \{\, x\in K \mid v(x-1)\geq t \,\} \qquad
      \text{(for $t\geq 0$).}
\]
Note that $\p_K^1=(p)$ and that $\p_K^{1/e}$ is the maximal ideal of $\OO_K$.

Our crucial assumption on $K$ is:

\begin{ass}\  \label{mainass}
\begin{itemize}
\item[(a)]
  $[\bar{K}:\bar{K}^p]=p$.
\item[(b)]
  $H^1(\Kb,\ZZ/p\ZZ)\neq 0$.
\end{itemize}
\end{ass}
Assumption \ref{mainass} is satisfied, for instance, if $\Kb$ is a function
field of transcendence degree one over a perfect field,
or if $\Kb$ is complete with respect to a discrete valuation and has a perfect
residue field (i.e.\ if $K$ is a {\em $2$-local field}).

\subsection{} \label{prelim2}

We let
\[
       F:=\bigcap_i \Kb^{p^i}
\]
denote the maximal perfect subfield of $\Kb$. Since $K$ is complete, the
inclusion $F\inj\Kb$ lifts to a unique embedding of the ring of Witt vectors
$W(F)$ into $\OO_K$, see \cite{SerreCL}, Proposition II.10 and Theorem
II.5. Let $k_0$ denote the fraction field of $W(F)$; we consider $k_0$ as a
subfield of $K$. Let $k$ denote the algebraic closure of $k_0$ in $K$. We call
$k$ the {\em field of constants} of $K$. Elements of $k$ are called {\em
  constants}.

\begin{prop} \label{constantprop}
\begin{enumerate}
\item The extension $k/k_0$ is finite and totally ramified. In particular, $k$
  is complete with respect to $v$ and has residue field $F$.
\item
  The field $k$ is the largest subfield of
  $K$ with perfect residue field.
\item
  Let $l/k$ be an algebraic extension. Then $l$ is the field of constants of
  $L:=Kl$.
\end{enumerate}
\end{prop}

\proof The residue field $\bar{k}$ of $k$ is contained in $\Kb$ and is an
algebraic extension of $F$. Since $F$ is the maximal perfect subfield of
$\Kb$, we have $F=\bar{k}$, i.e.\ the extension $k/k_0$ has trivial residue
field extension. Now it follows from \cite{SerreCL}, Theorem II.4 that $k/k_0$
is finite and totally ramified. So (i) is proved. 

Let $M\subset K$ be a subfield with perfect residue field $\Mb$. Let $M_0$
denote the fraction field of $W(\Mb)$.  Since $\Mb$ is a perfect subfield of
$\Kb$, we have $\Mb\subset F$ and hence $M_0\subset k_0$. But $M/M_0$ is
finite (use again \cite{SerreCL}, Theorem II.4). Therefore, $M\subset k$,
which proves (ii).

Let $l/k$ be an algebraic extension and set $L:=Kl$. We note that $l$ is
algebraically closed in $L$. Let $l'$ denote the field
of constants of $L$. Then $l\subset l'$, by (ii). The residue field
extension $\Lb/\Kb$ is algebraic. An easy argument shows that the maximal
perfect subfield $F':=\cap_i\Lb^{p^i}$ is algebraic over the maximal perfect
subfield $F$ of $\Kb$. But $F'$ is the residue field of $l'$ and $F$ is
contained in the residue field of $l$. It follows that $l'$ is
algebraic over $l$. We conclude $l=l'$. This completes the proof of the
proposition. 
\Endproof

\begin{defn} \label{constantdef} An extension $L/K$ is called {\em constant}
  if $L=K\, l$ for some algebraic extension $l/k$. An extension $L/K$ is
  called {\em almost constant} if it lies in the composition of a constant and
  an unramified extension.
\end{defn}

A second important assumption on $K$ that we will make throughout is that
the extension $K/k$ is {\em weakly unramified}, i.e.\ we have
\[
      v(k^\times)=v(K^\times).
\]
By a theorem of Epp (see \cite{Epp73}) this assumption always holds after
replacing $K$ by a finite constant extension. So in contrast to Assumption
\ref{mainass}, this is no real restriction of generality.

\subsection{}

An extension $L/K$ is called {\em fierce} if the extension of residue
fields $[\bar{L}:\bar{K}]$ is purely inseparable, and moreover
\[
        [\bar{L}:\bar{K}]=[L:K].
\]
Note that if $L/K$ is fierce and $K'/K$ is almost constant, then $L':=LK'/K'$
is again fierce, and we have $[L':K']=[L:K]$. Together with Epp's theorem,
this shows the following.

\begin{prop} \label{Eppprop}
  Let $L/K$ be a finite extension. Then there exists a finite, almost constant
  extension $K'/K$ such that the extension $L':=LK'/K'$ is
  fierce. Furthermore, the degree $[L':K']$ does not depend on the choice of
  $K'/K$. 
\end{prop}

In this paper we study the ramification of finite cyclic extensions $L/K$ `up
to almost constant extensions'. By the above proposition, it is no loss of
generality to assume from the start that $L/K$ is fierce.

\subsection{} \label{pit}

After replacing $K$ by a finite, constant extension, we may assume that $K$
contains a primitive $p$th root of unity $\zeta_p\in\mu_p(K)$. We set
$\lambda:=\zeta_p-1$. Note that $v(\lambda)=1/(p-1)$ and
\begin{equation} \label{lambdacongreq}
     \frac{\lambda^{p-1}}{p} \equiv -1 \pmod{\p_K^{1/e}}.
\end{equation}
It is clear that $\lambda$ is contained in the field of constants $k$ of
$K$. 

Let $\pi$ be a prime element of $k$ (which is also a prime element of
$K$!). Then $v(\pi)=1/e$, where $e$ is the absolute ramification index of
$k$. Note that $e$ is a multiple of $p-1$ since $\lambda\in k$. Let
$\cb\in\kb$ denote the residue class of the unit
$\pi^{e/(p-1)}/\lambda$. After replacing $k$ by a finite unramified extension
we may assume that there exists $\db\in\kb$ such that
$\db^{e/(p-1)}\cb=1$. Let $d\in\OO_k^\times$ be a lift of $\db$. We set
\[
     \pi_{1/e}:=d\cdot\pi.
\]
By definition $\pi_{1/e}$ is an element of $k$ with $v(\pi_{1/e})=1/e$ and 
\[
     \pi_{1/e}^{e/(p-1)} \equiv \lambda \pmod{\p_K^{1/(p-1)+1/e}}.
\]
Since the absolute ramification index of $k$ has not changed by passing to an
unramified extension, $\pi_{1/e}$ is a prime element of $k$. For any
$t\in v(K^\times)=e^{-1}\cdot\ZZ$ we set
\[
      \pi_t:=\pi_{1/e}^{et}.
\]
Now the following holds:
\begin{rem} \label{pitrem}
\begin{enumerate}
\item
  $v(\pi_t)=t$, for all $t\in v(K^\times)$.
\item
  $\pi_t^k=\pi_{k t}$, for all $k\in\ZZ$.
\item
  $\pi_{1/(p-1)} \equiv \lambda \pmod{\p_K^{1/(p-1)+1/e}}$.
\item
  $\pi_1 \equiv -p \pmod{\p_K^{1+1/e}}$.
\end{enumerate}
\end{rem}

\begin{rem} \label{pitrem2} 
  We will often have to replace $K$ by some finite constant extension
  $K'/K$. Then it may not be possible to extend the definition of $\pi_t$ to
  $K'$ in a compatible way. However, if we allow to further replace $K'$ by an
  unramified constant extension, then it is possible to define
  $\pi_t'$ for all $t\in v((K')^\times)$ in such a way that Remark
  \ref{pitrem} holds and 
  $\pi_t'/\pi_t\equiv 1$ for all $t\in v(K^\times)$. 

  As a result, if we replace $K$ by a constant extension then we need to
  choose $\pi_t$ once again, but this will never cause any problems.
\end{rem}

\section{Ramification and Swan conductors} \label{swan}

\subsection{} \label{swan1}

Let $L/K$ be a fierce Galois extension. Then the extension of valuation rings
$\OO_L/\OO_K$ is monogenic. In fact, the purely inseparable extension of
residue fields $\Lb/\Kb$ can be generated by one element $\xb\in\Lb$ by
Assumptions \ref{mainass}, and then any element $x\in\OO_L$ lifting $\xb$
generates $\OO_L$ over $\OO_K$. 

It follows that we have the usual theory of ramifications groups available,
see e.g.\ \cite{KatoVC}, \S 2-\S 3. Let $G=\Gal(L/K)$ denote the Galois group
of $L/K$. There are two filtrations by subgroups of $G$,
\[
     (G^t)_{t\geq 0}, \qquad (G_t)_{t\geq 0},
\]
with index set $\RR_{\geq 0}$ ({\em upper} and {\em lower numbering}), defined
as follows. 

For $\sigma\in G$, let
\[
       i_G(\sigma) := \min\,\{\,v(\sigma(x)-x) \mid x\in \OO_L\,\}.
\]
Then the lower numbering filtration is defined by
\[
      G_t := \{\,\sigma\in G \mid i_G(\sigma)\geq t \,\}.
\]
Note that this definition differs from the standard definition given e.g.\ in
\cite{SerreCL}, but it agrees with the definition used in \cite{KatoVC} (put
$\epsilon:=0$ in \cite{KatoVC}, Corollary 3.3). 

Set 
\[
     \phi_{L/K}(t) := \int_0^t\abs{G_s}\,ds.
\]
Then $\phi_{L/K}:\RR_{\geq 0}\to\RR_{\geq 0}$ is piecewise linear and strictly
increasing, hence bijective. Let $\psi_{L/K}:=\phi_{L/K}^{-1}$ be the inverse
map. The upper numbering filtration is now defined by
\[
   G^t:=G_{\psi_{L/K}(t)}.
\]
One shows that 
\[
     \psi_{L/K}(t) = \int_0^t \abs{G^s}^{-1}\,ds.
\]

\begin{rem} \label{swanrem}
The filtrations $(G^t)_t$ and $(G_t)_t$ have the following properties:
\begin{enumerate}
\item
  If $s\leq t$, then $G^s\supset G^t$ and $G_s\supset G_t$.
\item
  $G=G^0=G_0$.
\item
  For each $\sigma\in G\backslash\{1\}$, the sets $\{t\geq 0 \mid \sigma\in
  G^t\}$ and  $\{s\geq 0 \mid \sigma\in G_s\}$ have a maximum.
\item Let $H\lhd G$ be a normal subgroup and $M:=L^H$ the corresponding
  subfield. Then 
  \[
       (G/H)^t=G^t/(G^t\cap H), \qquad G^t\cap H = H^{\psi_{M/K}(t)},
  \]
  where the filtration $(H^t)_t$ is induced from the isomorphism
  $H\cong\Gal(L/M)$. 
\end{enumerate}
See e.g.\ \cite{KatoVC}, Lemma 2.9.
\end{rem}

It will be important for us that the ramification filtrations are invariant
under constant extensions of $K$.

\begin{lem} \label{ramlem1}
  Let $K'/K$ be an almost constant extension. Set $L':=LK'$ and
  $G':=\Gal(L'/K')$. Then the natural isomorphism $G'\iso G$ is compatible
  with the upper and lower filtrations, i.e.
  \[
      (G')_t\iso G_t, \quad (G')^t\iso G^t,
  \]
  for all $t\geq 0$.
\end{lem}

\proof
Let $x\in\OO_L$ be an element whose residue class generates the extension
$\Lb/\Kb$. Then $\OO_L=\OO_K[x]$ and $\OO_{L'}=\OO_{K'}[x]$. We conclude that 
\[
     i_G(\sigma) = v(\sigma(x)-x) = i_{G'}(\sigma),
\]
for all $\sigma\in G$.  The statement of the lemma follows immediately.
\Endproof

\subsection{} \label{swan2}

By a {\em
  character} on $K$ we shall always mean a continous group homomorphism
$\chi:\Gal(K^{\rm sep}/K)\to\QQ_p/\ZZ_p$. Let $L/K$ be the unique finite Galois
extension $L/K$ such that $\chi$ factors over $G=\Gal(L/K)$ and induces an
injective homomorpism $G\inj\QQ_p/\ZZ_p$. Then $G$ is cyclic of order $p^n$,
for some $n\geq 0$. We call $p^n$ the {\em order} of $\chi$. We call $\chi$
{\em weakly unramified} (resp.\ {\em fierce}) if the extension $L/K$ is weakly
unramified (resp.\ fierce). 

We are interested in ramification of characters, but only up to an almost
constant extension. By Proposition \ref{Eppprop} we may therefore assume
from the start that a given character on $K$ is fierce.

\begin{defn} \label{swandef}
  Let $\chi$ be a fierce character on $K$. We identify $\chi$ with a group
  homomorphism $\chi:G\to\QQ_p/\ZZ_p$, where $G=\Gal(L/K)$ is the Galois group
  of a suitable finite fierce Galois extension $L/K$ (but we do not assume
  that $\chi$ is injective). The (usual) {\em Swan conductor} of $\chi$ is the
  rational number
  \[
     \sw(\chi):=\max\,\{\,t\geq 0 \mid \chi|_{G^t}\neq 1\,\}.
  \]
  Note that $\sw(\chi)$ is independent of the chosen extension $L/K$ and does
  not change if we restrict $\chi$ to an almost constant extension of $K$. 
\end{defn}

Following \cite{KatoDSC}, we can also define a {\em refined Swan conductor}
of $\chi$, as follows. We define for $q,n\geq 1$ 
\[ 
     H_{p^n}^q(K):=H^q(K,\ZZ/p^n\ZZ(q-1))  
\]
and
\[
     H^q(K):=\varinjlim_n\, H_{p^n}^q(K).
\]
An element of $H^1(K)$ is simply a character on $K$ as defined
above, and 
\[
    H^2(K) = \Br(K)[p^\infty].
\]
Let 
\[
      (\cdot,\cdot)_K: H^1(K)\times K^\times \to H^2(K)
\]
be the pairing defined in \cite{SerreCL}, XIV, \S 1. We also need the two
morphism
\[
  i_1:\Kb\to H^1_p(K), \qquad i_2:\Omega_{\Kb}^1\to H^2_p(K)
\]
defined as follows. The morphism $i_1$ is the composition of the
Artin-Schreier map
\[
     \Kb\to H^1_p(\Kb):=H^1(\Kb,\ZZ/p\ZZ)
\]
with the restriction map
\[
      H^1_p(\Kb)\to H^1_p(K).
\]
Note that an element $\ab\in\Kb$ lies in the kernel of $i_1$ if and only if
$\ab=\bb^p-\bb$ for some $\bb\in\Kb$.

In order to define $i_2$ we choose an element $x\in\OO_K$ with
$\xb\not\in\Kb^p$. By Assumption \ref{mainass} (a), every differential
$\omega\in\Omega_{\Kb}^1$ can be written as
\[
    \omega= \ab\,\frac{d\xb}{\xb},
\]
for a unique element $\ab\in\Kb$. We set
\[
    i_2(\omega) := (i_1(\ab),x)_K.
\]
We can also write $\ab=\bb_0^p+\bb_1^p\xb+\ldots+\bb_{p-1}^p\xb^{p-1}$. The {\em
  Cartier operator} $\C:\Omega_{\Kb}\to\Omega_{\Kb}$ is the map defined by
\[
      \C(\omega) := \bb_0\,\frac{d\xb}{\xb},
\]
see \cite{CartierOp}.
It is shown in \cite{KatoGC} that $i_2(\omega)$ lies in $H^2_p(K)$ and does
not depend on the choice of $x$. Moreover, $i_2(\omega)=0$ if and only if
$\omega=\C(\eta)-\eta$ for some $\eta\in\Omega_{\Kb}^1$. 

 Using the above definition and Assumption \ref{mainass}
(b) one easily proves:

\begin{rem} \label{Crem}
\begin{enumerate}
\item
  The map
  \[
      \Omega_{\Kb}^1\to H^2_p(K), \qquad \omega \mapsto 
       \big(\;\zb\mapsto i_2(\zb\omega)\;\big)
  \]
  is $\Kb$-linear and injective.
\item
  We have
  \[
      i_2(\zb^p\omega)=i_2(\zb\,\C(\omega)),
  \]
  for all $\omega\in\Omega_{\Kb}^1$ and $\zb\in\Kb$.
\end{enumerate}
\end{rem}

\begin{thm}[Kato] \label{rswanthm} Let $\chi\in H^1(K)$ be a fierce character
  on $K$. Then the following holds.
  \begin{enumerate}
  \item
    The Swan conductor $\sw(\chi)$ is an element of
    $v(K^\times)=\ZZ\cdot\frac{1}{e}$.    
  \item
    For $s\geq 0$ the restriction of the morphism
    \[
        K^\times \to H^2(K), \qquad b\mapsto(\chi,b)_K
    \]
    to the subgroup 
    \[
        U_K^s := \{\,x\in K^\times \mid v(x-1)\geq s \,\}
    \] 
    is trivial if and only if $s>\sw(\chi)$.
  \item Let $\pi\in K$ be an element with $v(\pi)=\sw(\chi)$ (which exists
    by (i)). Then there exists a unique nonzero differential
    $\omega\in\Omega_{\Kb}^1$ such that
    \[
         (\chi,1-\pi z)_K = i_2(\zb\,\omega),
    \]
    for all $z\in\OO_K$.
  \end{enumerate}
\end{thm}

\proof See \cite{KatoDSC}, Theorem 3.6.
\Endproof

\begin{cor} \label{rswcor} 
  Let $\chi\in H^1(K)$ be a fierce
  character on $K$ of order $>1$. Then $\sw(\chi^p)<\sw(\chi)$.
\end{cor}

\begin{defn} \label{rswdef}
  Let $\chi\in H^1(K)$ be a fierce character on $K$, and set
  $\delta:=\sw(\chi)$. Let $\pi\in K$ be an element with $v(\pi)=\delta$ and 
  $\omega\in\Omega_{\Kb}^1$ as in Theorem \ref{rswanthm}. The {\em
    refined Swan conductor} of $\chi$ is the element
  \[
     \rsw(\chi):=\pi^{-1}\otimes \omega 
         \in \p_K^{-\delta}\otimes_{\OO_K}\Omega_{\Kb}^1
  \]
  (which does not depend on the choice of $\pi$).
\end{defn}

\begin{lem} \label{ramlem2}
  Let $K'/K$ be an almost constant extension. Then $\sw(\chi|_{K'})=\sw(\chi)$.
  Moreover, $\rsw(\chi|_{K'})$ is equal to the image of $\rsw(\chi)$ under
  the (injective!) morphism
  \[
      \p_K^{-\delta}\otimes_{\OO_K}\Omega_{\Kb}^1 \to
        \p_{K'}^{-\delta}\otimes_{\OO_{K'}}\Omega_{\Kb'}^1.
  \]
\end{lem}

\proof
The equality  $\sw(\chi|_{K'})=\sw(\chi)$ follows immediately from Lemma
\ref{ramlem1}. The second claim then follows from Theorem \ref{rswanthm} and
the formula
\[
    (\chi|_{K'},\,a)_{K'} = \Res_K^{K'} (\chi,\,a)_K, 
\]
for $a\in K^\times$ (see \cite{SerreCL}, XIV, \S 1).
\Endproof
 
The refined Swan conductor $\rsw(\chi)$, as it is defined above, is intrinsic
and encodes the usual Swan conductor $\sw(\chi)$. For our purposes, the
following definition, which depends on the choice of the elements $\pi_t$ in \S
\ref{pit}, will be more convenient.

\begin{defn} \label{dswdef}
  Let $\chi$ be a fierce character on $K$ of order $>1$ and set
  $\delta:=\sw(\chi)$. The {\em differential Swan conductor} of $\chi$ is the
  unique element $\omega=\dsw(\chi)\in\Omega_{\Kb}^1$ such that
  \[
      \rsw(\chi) = \pi_{\delta}^{-1}\otimes \omega
  \]
  (see \S \ref{pit} for the definition of $\pi_\delta$).  By definition, we
  have
  \begin{equation} \label{dsweq}
      (\chi,\,1-\pi_\delta\cdot z)_K = i_2(\zb\omega),
  \end{equation}
  for all $z\in\OO_K$. 
\end{defn}

\begin{rem}
  The definition of $\omega=\dsw(\chi)$ depends on the choice of $\pi_t$ made
  in \S \ref{pit}, but only up to a constant in $\kb^\times$. Moreover, the
  definition of $\omega=\dsw(\chi)$ is invariant under replacing $K$ by a
  finite constant extension, by Lemma \ref{ramlem2} and Remark \ref{pitrem2}.
\end{rem}

The following lemma shows how the refined Swan conductor behaves with respect
to addition of characters. 

\begin{lem} \label{rswlem}
  Let $\chi_i$, $i=1,2,3$, be fierce characters on $K$ satisfying the relation
  $\chi_3=\chi_1+\chi_2$. Set $\delta_i:=\sw(\chi_i)$ and
  $\omega_i:=\dsw(\chi_i)$, for $i=1,2,3$. Then the following holds.
  \begin{enumerate}
  \item
    If $\delta_1\neq\delta_2$ then
    $\delta_3=\max\{\delta_1,\delta_2\}$. Furthermore, we have
    $\omega_3=\omega_1$ if $\delta_1>\delta_2$ and $\omega_3=\omega_2$
    otherwise. 
  \item
    If $\delta_1=\delta_2$ and $\omega_1+\omega_2\neq 0$ then
    $\delta_1=\delta_2=\delta_3$ and $\omega_3=\omega_1+\omega_2$.
  \item
    If $\delta_1=\delta_2$ and $\omega_1+\omega_2=0$ then $\delta_3<\delta_1$.
  \end{enumerate}
\end{lem}

\proof
Follows from Theorem \ref{rswanthm} and Definition \ref{dswdef}. Details are
left to the reader.
\Endproof

\section{Cyclic extensions and ramification data}

Before we can state our main results, we need to define the notion of a {\em
  ramification datum}.

\begin{defn} \label{ramdatdef}
\begin{enumerate}
\item
  Let $n\geq 1$. 
  A {\em ramification datum} is a tupel
  $(\delta_i,\omega_i)_{i=1,\ldots,n}$, where $\delta_i\in\QQ_{\geq 0}$ is a
  nonnegative rational number and $\omega_i\in\Omega_{\Kb}^1\backslash\{0\}$
  is a nonzero differential form. 
\item
  Let $\chi\in H^1(K)$ be a fierce character on $K$ of order $p^n$, with
  $n\geq 1$. For $i=1,\ldots,n$ we set
  \[
       \chi_i:=p^{n-i}\cdot\chi, \qquad \delta_i:=\sw(\chi_i), 
         \qquad \omega_i:=\dsw(\chi_i).
  \]
  The tupel $(\delta_i,\omega_i)_{i=1,\ldots,n}$ is called the {\em
    ramification datum} associated to $\chi$.
  \end{enumerate}
\end{defn}

\begin{rem}
  Let $L/K$ be the cyclic extension of degree $p^n$ corresponding to the
  character $\chi$, and $G=\Gal(L/K)$. Let $(\delta_i,\omega_i)_i$ be the
  ramification datum associated to $\chi$. The following statements follow
  from Definition \ref{ramdatdef} and Corollary \ref{rswcor}. 
  \begin{enumerate}
  \item
    The numbers $\delta_i$ are precisely the breaks for the upper numbering
    filtration $(G^t)_t$. For $i=1,\ldots,n$  we have 
    \[
         \abs{G^{\delta_i}}=p^{n-i+1}.
    \]
    It follows that
    \[
        0<\delta_1<\delta_2<\ldots<\delta_n.
    \]
  \item
    It follows from (i) that the tupel $(\delta_i)$ only depends on the
    extension $L/K$. This is not quite true for the differentials
    $\omega_i$. If we replace $\chi$ by $\chi^a$, where $a\in\ZZ$ is prime to
    $p$, then $\omega_i$ gets replaced by $\bar{a}\,\omega_i$ (and where
    $\bar{a}$ denotes the residue of $a$ modulo $p$). 
  \end{enumerate}
\end{rem}

Here is our first main theorem.

\begin{thm} \label{thm1}
  Let $\chi\in H^1(K)$ be a fierce character of order $p^n>1$ and let
  $(\delta_i,\omega_i)_i$ denote the ramification datum associated to
  $\chi$. Then for all $i=1,..,n$ the following holds.
  \begin{enumerate}
  \item
    $0<\delta_1\leq p/(p-1)$. Moreover, we have
    \begin{enumerate}
    \item
      \quad $\delta_1=p/(p-1) \quad\Leftrightarrow\quad
                    \C(\omega_1)=\omega_1$,
    \item
      \quad $\delta_1<p/(p-1) \quad\Leftrightarrow\quad \C(\omega_1)=0$.
    \end{enumerate}
  \item
    Suppose $i>1$. If $\delta_{i-1}>1/(p-1)$ then we have
    \[
        \delta_i = \delta_{i-1}+1, \qquad
        \omega_i = -\omega_{i-1}.
    \]
  \item
    Suppose $i>1$. If $\delta_{i-1}\leq 1/(p-1)$ then
    \[
        p\,\delta_{i-1}\leq \delta_i \leq \frac{p}{p-1} .
    \]
    Moreover, we have
    \begin{enumerate}
    \item
      \quad $\delta_i=p\,\delta_{i-1}<p/(p-1) \quad\Rightarrow\quad
        \C(\omega_i)=\omega_{i-1}$,
    \item
      \quad $p\,\delta_{i-1}<\delta_i<p/(p-1) \quad\Rightarrow\quad
        \C(\omega_i)=0$,
    \item
      \quad $p\,\delta_{i-1}<\delta_i=p/(p-1) \quad\Rightarrow\quad
        \C(\omega_i)=\omega_i$,
    \item
      \quad $p\,\delta_{i-1}=\delta_i=p/(p-1)\quad\Rightarrow\quad
                 \C(\omega_i)=\omega_i+\omega_{i-1}$.
    \end{enumerate}
  \end{enumerate}
\end{thm}

\begin{rem} \label{thm1rem} In \cite{Hyodo87} Hyodo studies the ramification
  of cyclic Galois extensions $L/K$, where $K$ is a complete discretely valued
  field of mixed characteristic and imperfect residue field. His results imply
  some parts of Theorem \ref{thm1}. For simplicity we give details only in the
  case $n=2$ in Theorem. Namely, in the situation of Theorem \ref{thm1},
  \cite{Hyodo87}, Lemma 4.1, implies the following inequalities:
  \begin{itemize}
  \item[(a)]
    If $\delta_1\geq 1/(p-1)$ then
    \[
       \delta_1+1\leq \delta_2 \leq \frac{p}{p-1}+\delta_1\cdot\frac{p-1}{p}.
    \]
  \item[(b)]
    If $\delta_1\leq 1/(p-1)$ then
    \[
       p\,\delta_1\leq\delta_2\leq \frac{p}{p-1}+\delta_1\cdot\frac{p-1}{p}.
    \]
  \end{itemize}
  The statement of Theorem \ref{thm1} is stronger then (a) and (b)
  above. In fact, Theorem \ref{thm1} is false without Assumption
  \ref{mainass}, whereas Hyodo's results are valid in much greater
  generality.  
\end{rem}

\begin{rem}
  In the statement of Theorem \ref{thm1} (iii) the
  converse implications hold in Case (a) and (c), but not in Case
  (b) and (d). This can be seen as follows.

  Let us first consider (a) and assume that $\C(\omega_i)=\omega_{i-1}$.
  Since $\omega_j\neq 0$ for all $j$ by definition, we conclude that
  $\C(\omega_i)\neq 0,\omega_i+\omega_{i-1}$. This rules out Case (b) and
  (d). So suppose that (c) holds and therefore
  $\omega_i=\C(\omega_i)=\omega_{i-1}$. Using $\delta_j\leq
  p^{-1}\delta_i<p/(p-1)$ for $j<i$ and the forward direction of (a) one shows
  inductively that $\omega_i=\omega_{i-1}=\ldots=\omega_1$. But then
  $\delta_1<p/(p-1)$ and $\C(\omega_1)=\omega_1$. This contradicts Part (i) of
  Theorem \ref{thm1} and proves that the converse implication holds in Case
  (a). The proof in Case (c) is similar.

  On the other hand, it is possible that
  $\C(\omega_i)=\omega_i+\omega_{i-1}=0$ which shows that the converse
  implication can't hold in Case (b) and (d).
\end{rem}

Our second main result states that every ramification datum satisfying the
conclusion of Theorem \ref{thm1} is realized by some fierce
character. More precisely:

\begin{thm} \label{thm2} Let $(\delta_i,\omega_i)_i$ be a ramification datum
  satisfying Conditions (i)-(iii) in Theorem \ref{thm1}. Then there exists a
  finite constant extension $K'/K$ and a fierce character
  $\chi$ of order $p^n$ on $K'$, such that $(\delta_i,\omega_i)$ is the
  ramification datum associated to $\chi$.
\end{thm} 

The proof of Theorem \ref{thm1} and Theorem \ref{thm2} occupies the rest of
this paper. Both theorems are proved by induction on $n$. The case $n=1$ is
proved in Section \ref{ordp}, and most of the induction step is done in Section
\ref{ind}. Finally, a crucial case occuring in the induction step for Theorem
\ref{thm2} is dealt with in Section \ref{minimal}.

\section{Cyclic extensions of order $p$} \label{ordp}

In this section we prove Theorem \ref{thm1} and Theorem \ref{thm2} for fierce
characters of order $p$. The proof relies on an explicit description of the
refined Swan conductor $\rsw(\chi)$ of such characters in terms of Kummer
theory. This description will also be useful for the study of characters of
higher order.

Kummer theory defines an isomorphism
\[
    K^\times/(K^\times)^p \iso H_p^1(K).
\]
We write $\chi=\chi_u\in H^1_p(K)$ for the character corresponding to an element
$u\in K^\times$. Explicitly, $\chi$ is given as follows. If $u\in K^p$, then
$\chi$ is the trivial character. Otherwise, let $v\in K^{\rm alg}$ be a
$p$th root of $u$ and $L:=K[v]$. Clearly, $L/K$ is a Galois extension, and the
map
\[
      G=\Gal(L/K)\to \mu_p(K), \quad \sigma\mapsto \sigma(v)/v,
\]
is an isomorphism of cyclic groups of order $p$. The value
$\chi(\sigma)\in\ZZ/p\ZZ$ is determined by the identity
\[
     \sigma(v)/v = \zeta_p^{\,\chi(\sigma)p}.
\]
Note that $\chi$ depends on the choice of $\zeta_p$ but not on the choice of
$v$. 

\begin{defn} \label{reduceddef}
  An element $u\in K^\times$ is said to be {\em reduced} if one of the
  following two conditions hold.
  \begin{itemize}
  \item[(a)]
    $u\in\OO_K^\times$ and $\ub\not\in\Kb^p$.
  \item[(b)]
    $u=1+\pi^pw$, with $\pi\in k$, $0<v(\pi)<1/(p-1)$, $w\in\OO_K$ and
    $\wb\not\in\Kb^p$. 
  \end{itemize}
  We say that $u$ is {\em reducible} if there exists $a\in K^\times$ such that
  $ua^p$ is reduced. 
\end{defn}
     
\begin{prop} \label{reducedprop}
  The character $\chi=\chi_u$ is fierce if and only if $u$ is reducible.
\end{prop}

\proof This follows from \cite{Hyodo87}, Lemma 2.16. For convenience of the
reader we sketch a proof of the `if'-direction. By definition, the character
$\chi$ corresponds to the Galois extension $L:=K[v]/K$, were $v^p=u$. Assume
that $u$ is reducible. In order to show that $L/K$ is fierce we may even
assume that $u$ is reduced.  In Case (a), the residue $\vb\in\Lb$ generates an
inseparable extension of $\Kb$ of degree $p$. It follows that
$[L:K]=[\Lb:\Kb]_{\rm ins}$, i.e.\ $L/K$ is fierce. Note also that
$\OO_L=\OO_K[v]$ in this case.

In Case (b) we set $z:=(v-1)/\pi\in L$. Then
\begin{equation} \label{reducedpropeq1}
    \frac{(\pi z+1)^p-1}{\pi^p} = z^p+p\pi^{-1}z^{p-1}+\ldots+p\pi^{1-p} z = w.
\end{equation}
By the assumption on $\pi$, the valuation of the coefficient of $z^i$ in
\eqref{reducedpropeq1} is
\[
   v\big(\,\binom{p}{i}\pi^{i-p} \,\big)
     = 1-(p-i)v(\pi) > 0, \quad\text{for $i=1,\ldots,p-1$.}
\]   
It follows that $z\in\OO_L$ is integral and that its residue satifies the
irreducible equation 
\[
      \zb^p=\wb.
\]
As in Case (a), we conclude that $L/K$ is fierce and that $\OO_L=\OO_K[z]$.
\Endproof

\begin{prop} \label{reducedprop2}
  Let $\chi=\chi_u\in H^1_p(K)$, where $u\in K^\times$ is reduced. Then the
  refined Swan conductor $\rsw(\chi)$ is given as follows.
  \begin{enumerate}
  \item
    If $\ub\not\in\Kb^p$ (Case (a) in Definition \ref{reduceddef}) then 
    \[
          \rsw(\chi) = \lambda^{-p}\otimes \frac{d\ub}{\ub}.
    \]
  \item
    If $u=1+\pi^pw$ (Case (b) in Definition \ref{reduceddef}) then
    \[
        \rsw(\chi) = (\pi^p\lambda^{-p})\otimes d\wb.
    \]
  \end{enumerate}
\end{prop}

\proof We use the notation of the proof of Proposition \ref{reducedprop}. Let
$\sigma\in G=\Gal(L/K)$ be the element with $\chi(\sigma)=1/p$, i.e.\
$\sigma(v)=\zeta_p\, v$. Fix an element $x\in \OO_L^\times$ with
$\Lb=\Kb[\xb]$. Write $y:=N_{L/K}(x)\in\OO_K$. Then $\yb=\xb^p\in\Kb\backslash
\Kb^p$. Set $a:=\sigma(x)/x-1$, $b:=N_{L/K}(a)$. By \cite{KatoLCFT2}, \S 3.3,
Lemma 15, we have
\begin{equation} \label{reducedpropeq2}
  (\chi,\ 1-bc)_K =  (i_1(\cb),\, y)_K 
        = i_2\big(\,\cb\frac{d\yb}{\yb}\,\big),
\end{equation}
for all $c\in\OO_K$. It follows that
\begin{equation} \label{reducedpropeq3}
  \rsw(\chi) = b^{-1}\otimes \frac{d\yb}{\yb}.
\end{equation}

In Case (a), we set $x:=v$. Then $a=\zeta_p-1=\lambda$, $b=\lambda^p$ and
$y=u$. Hence $\rsw(\chi)=\lambda^{-p}\otimes d\ub/\ub$ by
\eqref{reducedpropeq3}, and (i) is proved.

In Case (b), we set $x:=z$. Then $a=\lambda\pi^{-1}z^{-1}+\lambda$ and $y=w$
(if $p=2$ then $y=-w$). The assumption on $\pi$ implies
$0<v(a)=1/(p-1)-v(\pi)<v(\lambda)$ and
\[
    a \equiv \lambda\pi^{-1}z^{-1} \pmod{p^{v(a)+\epsilon}}.
\]
It follows that
\[
    b \equiv \lambda^p\pi^{-p}w^{-1} \pmod{p^{pv(a)+\epsilon}}.
\]
(Note that the case $p=2$ is no exception since $-1=1$ in $\Kb$.)  
As in Case (a) we conclude that
\[
  \rsw(\chi) = b^{-1}\otimes d\wb/\wb = (\pi^p\lambda^{-p})\otimes d\wb.
\]
\Endproof

\begin{cor} \label{n1cor}
\begin{enumerate}
\item
  Let $(\delta,\omega)$ be the ramification datum associated to a fierce
  character $\chi$ of order $p$. 
  Then we have either 
  \[
       \delta=\frac{p}{p-1},\qquad \omega=\frac{d\yb}{\yb}, 
  \]
  or 
  \[
      0<\delta<\frac{p}{p-1}, \qquad \omega = d\yb,
  \]
  for some element $\yb\in\Kb\backslash\Kb^p$.
\item
  Theorem \ref{thm1} and Theorem \ref{thm2} hold true for $n=1$.
\end{enumerate}
\end{cor}

\proof Claim (i) follows immediately from \eqref{lambdacongreq} and
Proposition \ref{reducedprop2}. For the proof of (ii), recall the well known
fact that for a differential form $\omega\in\Omega_{\Kb}^1$ we have
$\C(\omega)=\omega$ if and only if $\omega=d\yb/\yb$, for some element
$\yb\in\Kb$. Similarly, $\C(\omega)=0$ if and only if $\omega=d\yb$. In both
cases, $\omega\neq 0$ if and only if $\yb\not\in\Kb^p$.  It is now immediate
that Theorem \ref{thm1} follows from (i) for $n=1$.

Conversely, let $(\delta,\omega)$ be a deformation datum of length $n=1$ which
satisfies Condition (i) of Theorem \ref{thm1}. In Case (a), we can write
$\omega=d\ub/\ub$ for some $\ub\in\Kb\backslash\Kb^p$. Then for any lift
$u\in\OO_K$ of $\ub$, the deformation datum associated to the character
$\chi_u$ is equal to $(\delta,\omega)$. Case (b) is similar, the only
difference being that we need that the rational number $\delta$ is contained
in $v(K^\times)$, which holds after a finite constant extension of $K$. Claim
(ii) is now proved.
\Endproof 

\section{Higher order} \label{ind}

The case $n=1$ of Theorem \ref{thm1} and Theorem \ref{thm2} having been
proved, we continue by induction with the case $n>1$. 

Theorem \ref{thm1} follows immediately, by induction, from the case $n=1$ and
the following proposition.

\begin{prop} \label{indprop1}
  Let $\chi$ be a fierce character of order $p^n$, $n>1$. Set
  $\chib:=p\cdot\chi$, 
  \[
      \delta:=\sw(\chi), \quad \omega:=\dsw(\chi), \quad
      \deltab:=\sw(\chib), \quad \omegab:=\dsw(\chib).
  \]
  Then the following holds.
  \begin{enumerate}
  \item
    $\delta\leq \max\{\,\deltab+1,p/(p-1)\,\}$.
  \item
    Suppose $\deltab>1/(p-1)$. Then
    \[
        \delta=\deltab+1, \qquad \omega=-\omegab.
    \]
  \item
    Suppose $\deltab\leq 1/(p-1)$. Then
    \[
         p\,\deltab \leq \delta \leq \frac{p}{p-1}.
    \]
    Moreover,
    \begin{enumerate}
    \item
      $\delta=p\,\deltab<p/(p-1) \quad\Rightarrow\quad
          \C(\omega)=\omegab.$
    \item
      $p\,\deltab<\delta<p/(p-1) \quad\Rightarrow\quad
          \C(\omega)=0.$
    \item
      $p\,\deltab<\delta=p/(p-1) \quad\Rightarrow\quad
          \C(\omega)=\omega.$
    \item
      $p\,\deltab=\delta=p/(p-1) \quad\Rightarrow\quad
          \C(\omega) = \omega+\omegab$.
    \end{enumerate}
  \end{enumerate}
\end{prop}

\proof
The relation $p\cdot\chi=\chib$ and the bilinearity of the symbol
$(\cdot,\cdot)_K$ imply the relation
\begin{equation} \label{indpropeq1}
     (\chib,a)_K = (\chi,a^p)_K,
\end{equation}
for all $a\in K^\times$. For an element of the form $a=1-\pi_tz$, with
$t\in v(\p_K)$ and $z\in\OO_K$, we get, using Remark \ref{pitrem}: 
\begin{equation} \label{indpropeq2}
    a^p = \;\begin{cases}
        1 - \pi_{pt}\cdot z^p + \ldots ,& \quad t<1/(p-1),\\
        1 - \pi_{pt}\cdot(z^p-z) + \ldots, & \quad t=1/(p-1),\\
        1 + \pi_{t+1}\cdot z + \ldots, & \quad t>1/(p-1).
            \end{cases}
\end{equation}
Here the dots indicate terms of higher valuation than the preceeding term. 

In order to prove (i), we assume that $\delta>\deltab+1$ and
$\delta>p/(p-1)$. Then $\delta-1>\deltab$, so Theorem \ref{rswanthm} (ii)
shows that
\[
     (\chib,\,1-\pi_{\delta-1}\cdot z)_K =0,
\]
for all $z\in\OO_K$. Using \eqref{indpropeq1}, \eqref{indpropeq2}, the
inequality $\delta-1>1/(p-1)$ and Theorem
\ref{rswanthm} (ii)-(iii) we deduce
\[\begin{split}
   0 & = (\chi,\,(1-\pi_{\delta-1}\cdot z)^p)_K \\
     & = (\chi,\, 1+\pi_\delta\cdot z + \ldots )_K \\
     & = (\chi,\, 1+\pi_\delta\cdot z)_K,
\end{split}\] 
for all $z\in\OO_K$. 
But this is a contradiction to Theorem \ref{rswanthm} (ii). We have proved (i). 

We proceed to prove (ii). By (i) and the assumption we have
$\delta\leq\deltab+1$. By the definition of $\deltab$ and $\omegab$ we have
\[
   (\chib,\,1-\pi_{\deltab}\cdot z)_K = i_2(\zb\omegab)\neq 0,
\]
for some $z\in\OO_K$. Using \eqref{indpropeq1} and \eqref{indpropeq2} we
deduce
\begin{equation} \label{inpropeq3} 
     (\chi,\,(1-\pi_{\deltab}\cdot z)^p)_K = 
      (\chi, 1+\pi_{\deltab+1}\cdot z +\ldots)_K\neq 0.
\end{equation}
Now the definition of $\delta$ and Theorem \ref{rswanthm} (ii) shows that
$\delta\geq\deltab+1$. Hence $\delta=\deltab+1$, and the above calculation shows
that 
\[
   i_2(\zb\omegab)=(\chi,1+\pi_\delta\cdot z)_K = i_2(-\zb\omega),
\]
for all $z\in\Kb$. Hence $\omega=-\omegab$, and (ii) is proved. 

The proof of (iii) follows the same line of argument. First of all, (i) and
the assumption imply $\delta\leq p/(p-1)$. We have
\begin{equation} \label{indpropeq4}
\begin{split}
   i_2(\zb\omegab) & = (\chi,\,(1-\pi_{\deltab}\cdot z)^p)_K  \\ 
     & = \begin{cases}
           \;(\chi,\,1-\pi_{p\deltab}\cdot z^p+\ldots)_K, & 
                                             \quad \deltab<1/(p-1),\\
           \;(\chi,\,1-\pi_{p\deltab}\cdot (z^p-z)+\ldots)_K, & 
                                             \quad \deltab=1/(p-1).
        \end{cases}
\end{split}
\end{equation}
Since $i_2(\zb\omegab)\neq 0$ for some $z$, we conclude that $\delta\geq
p\,\deltab$. 

Suppose that $\delta=p\,\deltab<p/(p-1)$. Then from \eqref{indpropeq4} and
Remark \ref{Crem} we get
\[
     i_2(\zb\omegab)=i_2(\zb^p\omega)=i_2(\zb\,\C(\omega)),
\]
for all $\zb\in\Kb$. This shows that $\C(\omega)=\omegab$ (Case
(a)). Similarly, if $p\,\deltab=\delta=p/(p-1)$, we get
\[
   i_2(\zb\omegab) = i_2((\zb^p-\zb)\omega)=i_2(\zb(\C(\omega)-\omega)),
\]
which shows $\C(\omega)=\omega+\omegab$ (Case (d)). 

Finally, suppose $p\,\deltab<\delta\leq p/(p-1)$. To prove the remaining cases
(b) and (c), we may assume that $\delta/p\in v(K^\times)$. Then
\[\begin{split}
    0 & = (\chib,\,1-\pi_{\delta/p}\cdot z)_K \\
      & = (\chi,\, (1-\pi_{\delta/p}\cdot z)^p)_K \\
      & = (\chi,\, 1-\pi_\delta\cdot g(z)+\ldots)_K = i_2(g(\zb)\,\omega), 
\end{split}\]
with $g(z):=z^p$ if $\delta<p/(p-1)$ and $g(z)=z^p-z$ if $\delta=p/(p-1)$. 
In the first case (Case (b)) we get
\[
     i_2(\zb\,\C(\omega))=i_2(\zb^p\omega)=0
\]
and conclude $\C(\omega)=0$. In the second case (Case (c)), we get
\[
     i_2(\zb\,(\C(\omega)-\omega))=i_2((\zb^p-\zb)\,\omega)=0
\]
and we conclude $\C(\omega)=\omega$. Now the proposition and Theorem
\ref{thm1} are proved. 
\Endproof

We now start the proof of Theorem \ref{thm2}. We fix a character $\chib$ on
$K$ of order $p^{n-1}$ and we assume that $\chib$ is fierce. We set
\[
     \deltab:=\sw(\chib), \qquad \omegab:=\dsw(\chib).
\]

\begin{defn} \label{liftdef}
By a {\em lift} of $\chib$ we mean a character $\chi$, defined (and fierce)
over some finite constant extension $K'/K$, with $p\cdot\chi=\chib|_{K'}$. 
\end{defn}

Theorem
\ref{thm2} follows by induction from the case $n=1$ and the following
proposition.

\begin{prop} \label{indprop2}
  Given $\delta\in\QQ_{>0}$ and $\omega\in\Omega_{\Kb}^1\backslash\{0\}$
  satisfying Condition (i)-(iii) of Proposition \ref{indprop1}, there exist a
  lift  $\chi$ of $\chib$ with $\delta=\sw(\chi)$ and $\omega=\dsw(\chi)$.
\end{prop}

\proof
After a finite constant extension we may assume that $K$ contains a $p^n$th
root of unity. Then Kummer theory shows that the following sequence is exact:
\[
     0 \to H^1_p(K)\longrightarrow H^1_{p^n}(K) 
            \lpfeil{\cdot p} H^1_{p^{n-1}}(K) \to 0.
\]
Therefore, the set of all lifts $\chi$ of $\chib$ (defined over $K$) is a
principal homogenous spaces under the natural action of $H^1_p(K)$. In
particular, it is nonempty. Note that, given an individual lift $\chi$, we can
make it fierce by a finite constant extension of $K$. However, we can't do
that for all lifts $\chi$ at a time.
 
Suppose first that $\deltab>1/(p-1)$. Then for any (fierce) lift $\chi$ we
have $\sw(\chi)=\deltab+1$ and $\dsw(\chi)=-\omegab$ by Proposition
\ref{indprop1} (ii). Since at least one lift exist, which becomes fierce after
a finite constant extension of $K$, there is nothing more to show.

Next suppose that $\deltab=1/(p-1)$, and let $\chi_0$ be some fierce lift of
$\chib$. Set $\delta_0:=\sw(\chi_0)$ and $\omega_0:=\dsw(\chi_0)$. Then we
have $\delta_0=p/(p-1)$ and $\C(\omega_0)=\omega_0+\omegab$ by Proposition
\ref{indprop1} (iii.d). Our assumption on $\delta$ and $\omega$ mean that
$\delta=p/(p-1)$ and $\C(\omega)=\omega+\omegab$. Set
$\eta:=\omega-\omega_0$. Then $\C(\eta)=\eta$. By the Case $n=1$ of Theorem
\ref{thm1} (see Corollary \ref{n1cor} (ii)) there exists a character $\psi$ of
order $p$ with $\sw(\psi)=p/(p-1)$ and $\dsw(\psi)=\eta$. Set
$\chi:=\chi_0+\psi$. This is a character of order $p^n$. After a constant
extension of $K$ it is fierce, and then Lemma \ref{rswlem} (ii) shows that
\[
    \sw(\chi)=\delta, \qquad \dsw(\chi)=\omega_0+\eta=\omega.
\]
So in this case the proposition is proved, too.

Finally we deal with the case $\deltab<1/(p-1)$. The hard part is to show the
following

\begin{lem} \label{minlem} Let $\chib\in H^1_{p^{n-1}}(K)$ be a fierce
  character of order $p^{n-1}$, with $\deltab=\sw(\chib)<1/(p-1)$. Then there
  exists a lift $\chi^{\rm min}$ of $\chib$ (defined over a finite constant
  extension $K'/K$) such that $\sw(\chi^{\rm min})=p\deltab$.
\end{lem}

The proof of Lemma \ref{minlem} will be given in the next section. Let
$\chi^{\rm min}$ be a lift of $\chib$ with $\sw(\chi^{\rm min})=p\deltab$. Set
$\omega^{\rm min}:=\dsw(\chi^{\rm min})$. We have $\C(\omega^{\rm
  min})=\omegab$ by Proposition \ref{indprop1} (iii.a). If $\delta=p\deltab$,
then we set $\eta:=\omega-\omega^{\rm min}$. Since
$\C(\eta)=\omegab-\omegab=0$ and $\delta<p/(p-1)$, there exists a constant
extension $K'/K$ and $\psi\in H^1_p(K')$ with $\sw(\psi)=\delta$ and
$\dsw(\psi)=\eta$ (Corollary \ref{n1cor} (ii)). Set $\chi:=\chi^{\rm
  min}+\psi$. It follows from Lemma \ref{rswlem} that $\sw(\chi)=\delta$ and
$\dsw(\chi)=\omega$, as required. 

The case $p\deltab<\delta$ is handled in a similar way. We are in Case (a) or
(c) of Proposition \ref{indprop1}. In Case (a) we have $\delta<p/(p-1)$ and
$\C(\omega)=0$, in Case (c) we have $\delta=p/(p-1)$ and
$\C(\omega)=\omega$. In both cases, Corollary \ref{n1cor} (ii) shows that
there exists a constant extension $K'/K$ and $\psi\in H^1_p(K')$ with
$\sw(\psi)=\delta$ and $\dsw(\psi)=\omega$. Set $\chi:=\chi^{\rm
  min}+\psi$. Using again Lemma \ref{rswlem} we see that $\sw(\chi)=\delta$
and $\dsw(\chi)=\omega$. This completes the proof of Proposition
\ref{indprop2}, under the condition that Lemma \ref{minlem} holds.
\Endproof

\section{Construction of a minimal lift} \label{minimal}

The goal of this section is to give a proof of Lemma \ref{minlem} and thus
complete the proof of Theorem \ref{thm2}.  Since the proof is quite
technical, we start by explaning the main idea.

Let $\chib$ be a character of order $p^{n-1}$, with
$\deltab:=\sw(\chib)<1/(p-1)$. 
Let $\chi_0$ be an arbitrary lift of
$\chib$. Then $\delta_0:=\sw(\chi_0)\geq p\deltab$, by Proposition
\ref{indprop1}. If $\delta_0=p\deltab$ then we are done. Otherwise, we can use
Lemma \ref{rswlem} to find a lift $\chi_1$ of $\chib$ with
$\delta_1:=\sw(\chi_1)<\delta_0$. This argument shows the following. If the
set 
\[
   \Delta := \{\,\sw(\chi) \mid \text{$\chi$ is a lift of $\chib$}\,\}
\]
has a minimum, then this minimum is equal to $p\deltab$, and we are
done. However, it is not obvious that a minimum exists: since lifts are
defined over finite but arbitrarily large constant extensions of $K$, $\Delta$
is not contained in any discrete subgroup of $\RR$.

Our solution to this problem is to consider a certain subset of the set of all
lifts of $\chib$. Let us call lifts that lie in this subset {\em moderate} (in
the actual proof this terminology is used in a slightly different way, see
Definition \ref{moderatedef}). We then show that the corresponding subset
$\Delta^{\rm mod}\subset\Delta$ has the following two properties: (a) if
$\Delta^{\rm mod}$ has a minimum, then it is equal to $p\deltab$, (b) if
$\delta_0>\delta_1>\ldots$ is a strictly decreasing sequence in $\Delta^{\rm
  mod}$, then $\lim_i \delta_i<p\deltab$ (again, in the actual proof this is
formulated differently). Combining (a) and (b) shows that there exists a minimal
lift $\chi$ with $\sw(\chi)=p\deltab$, proving Lemma \ref{minlem}.

\subsection{} \label{epsilonsubsec}

We fix a fierce character $\chib\in H^1_{p^{n-1}}(K)$ of order $p^{n-1}$, with
$n\geq 2$. It gives rise to a cyclic Galois extension $M/K$, with Galois group
$\bar{G}:=\Gal(M/K)\cong\ZZ/p^{n-1}$. Let $\delta_1,\ldots,\delta_{n-1}$ be
the breaks of the upper numbering filtration $(\bar{G}^t)_t$. Then
\[
    |\bar{G}^t| =\; \begin{cases}
              \;\;p^{n-1}, & \qquad t\leq \delta_1, \\
              \;\;p^{n-i}, & \qquad \delta_{i-1}<t\leq \delta_i,\\
              \;\;1      , & \qquad t>\delta_{n-1}.
            \end{cases}
\]
By definition, the Swan conductor $\deltab:=\sw(\chib)$ of $\chib$ is equal to
$\delta_{n-1}$. We assume that $\deltab<1/(p-1)$.

Let $\psi_{M/K}:\RR_{\geq 0}\to\RR_{\geq 0}$ be the inverse
Herbrand function associated to $M/K$, see \S \ref{swan1}. Then for
$t\geq\deltab$ we have
\begin{equation}  \label{psieq}
\begin{split}
   \psi_{M/K}(t) & = \delta_1\,p^{-n+1}+(\delta_2-\delta_1)\,p^{-n+2}+\ldots
                       +(\delta_{n-1}-\delta_{n-2})\,p^{-1} + (t-\deltab) \\
          & = t - \epsilon,
\end{split}
\end{equation}
with
\begin{equation} \label{epsiloneq}
   \epsilon := \delta_1\cdot\frac{p-1}{p^{n-1}}+
      \delta_2\cdot\frac{p-1}{p^{n-2}}+\ldots+\delta_{n-1}\cdot\frac{p-1}{p}>0.
\end{equation}
We may assume that $\epsilon\in v(K^\times)$. It is obvious from
\eqref{psieq} that
\begin{equation} \label{epsilonineq}
  \epsilon \leq \deltab<1/(p-1)\leq 1.
\end{equation}
It follows that 
\[
    \tilde{\OO}_M^{(\epsilon)}:=\OO_M^p+\pi_\epsilon\cdot\OO_M 
\]
is a subring of $\OO_M$.

\begin{lem} \label{epsilonlem} Let $a\in\OO_K$ be given. After replacing $K$
  by a finite constant extension (which may depend on $a$), we have
  $a\in\tilde{\OO}_M^{(\epsilon)}$.
\end{lem}

\proof Set $K':=K(a^{1/p})$, $M'=MK'=M[a^{1/p}]$. We may assume that, after a
finite constant extension of $K$, $M'/M$ is fierce of degree $p$ (otherwise, $a$
becomes a $p$th power in $\OO_M$ and we are done). But then $M'/K$ is a fierce
Galois extension such that
\[
     \Gal(M'/K) = \Gal(M'/K')\times \Gal(M'/M) \cong 
           \ZZ/p^{n-1}\ZZ\times\ZZ/p\ZZ.
\]
Let $\delta_{K'/K}$ (resp.\ $\delta_{M'/M}$) denote the Swan conductor of a
character of order $p$ giving rise to the extension $K'/K$ (resp.\ $M'/M$). 
Using Remark \ref{swanrem} and Definition \ref{swandef} it is easy to see 
that
\begin{equation} \label{epsilonlemeq1}
    \delta_{M'/M} = \psi_{M/K}(\delta_{K'/K}).
\end{equation}
By Proposition \ref{reducedprop} we can write 
\[
        a=b^p+\pi_t\cdot c,
\]
with $b,c\in\OO_K$, $\bar{c}\not\in\Kb^p$ and $t\in v(\OO_K)$. Then 
\[
       \delta_{K'/K} = p/(p-1)-t
\] 
by Proposition \ref{reducedprop2}. If $t\geq\epsilon$ then we are done. We may
therefore assume 
\[
     \delta_{K'/K} > p/(p-1)-\epsilon > 1/(p-1).
\]
But then \eqref{psieq} and \eqref{epsilonlemeq1} show that
\[
      \delta_{M'/M} = \delta_{K'/K}-\epsilon = p/(p-1)-t',
\]
with $t':=t+\epsilon\geq \epsilon$. Using again Proposition \ref{reducedprop}
and Proposition \ref{reducedprop2}, we see that we can write
\[
      a=(b')^p+\pi_{t'}\cdot c',
\]
with $b',c'\in\OO_M$. This proves the lemma.
\Endproof

\begin{cor} \label{epsiloncor}
  After replacing $K$ by a finite constant extension, the following holds. For
  all $\bar{a}\in\Kb$ there exists a lift $a\in\OO_K$ which can be written in
  the form
  \[
       a=b^p+\pi_\epsilon\cdot c,
  \]
  with $b,c\in\OO_M$.
\end{cor}

\proof
Choose $x\in\OO_K$ with $\xb\not\in\Kb^p$. Then 
\[
      \Kb=\Kb^p[\xb],
\]
by Assumption \ref{mainass} (a).  By Lemma \ref{epsilonlem} we may assume that
$x\in\tilde{\OO}_M^{(\epsilon)}$. It follows that $\Kb=\Kb^p[\xb]$ is
contained in the image of the ring homomorphism
\[
     \tilde{\OO}_M^{(\epsilon)}\cap\OO_K\to\Mb,
\]
proving the claim.
\Endproof

After replacing $K$ by a finite constant extension once, we may and will
assume from now on that the conclusion of Corollary \ref{epsiloncor} holds. 

\subsection{}

Let $\Lambda:=v(K^\times)$ denote the value group of $K$. Then
$\Lambda=\ZZ\cdot\frac{1}{e}$, where $e$ is the absolute ramification index of
$K$. For $k\geq 1$ we set
\[
     \nu_k:=1+\frac{1}{p}+\ldots+\frac{1}{p^{k-1}}.
\]
Note that $\nu_k\to p/(p-1)$ for $k\to\infty$.

\begin{defn} \label{Ledef}
  The subset $\Lambda_\epsilon\subset\QQ$ is defined as follows. A rational
  number $t\in\QQ$ is contained in $\Lambda_\epsilon$ if and only if $t\geq 0$
  and the implication
  \[
     t < \big(1-\epsilon/p\big)\cdot\nu_k \quad\Rightarrow\quad
           p^k\,t\in \Lambda 
  \]
  holds, for all $k\geq 1$. Here $\epsilon\in \Lambda$ is defined by
  \eqref{epsiloneq}. 
\end{defn}

It is clear that $\Lambda\cap\QQ_{\geq
  0}\subset\Lambda_\epsilon$. Furthermore:

\begin{lem} \label{Lelem}
  \begin{enumerate}
  \item
    $(\Lambda_\epsilon,+)$ is a submonoid of $(\QQ,+)$. 
  \item Suppose $t\in\Lambda_\epsilon$ and $t\geq\epsilon$. Then
    $t'_i:=1+i(t-\epsilon)/p$ and $t_i:=1+it/p$, for $i=1,2,\ldots$, lie in
    $\Lambda_\epsilon$ as well.
  \item Let $t_1,t_2,\ldots\in\Lambda_\epsilon$ be a strictly increasing
    sequence. Then
    \[
        \limsup_i t_i \geq \frac{p-\epsilon}{p-1}.
    \]
  \item
    Suppose $s,t\in\Lambda_\epsilon$, with $s,t\geq \epsilon$. Then
    $s+t-\epsilon\in\Lambda_\epsilon$. 
  \end{enumerate}
\end{lem}

\proof
Part (i) of the lemma follows directly from the definition of
$\Lambda_\epsilon$. To prove (ii) we first note that
\begin{equation} \label{Lelemeq1}
  p\cdot\big((1-\frac{\epsilon}{p})\nu_k -1\big) +\epsilon = 
    (1-\frac{\epsilon}{p})\nu_{k-1}.
\end{equation}
Assume that 
\[
      t_i'=1+\frac{i(t-\epsilon)}{p}<(1-\frac{\epsilon}{p})\nu_k.
\]
Since $t_1'\leq t_i'$ this inequality holds in particular for $i=1$.  
Using \eqref{Lelemeq1} and a direct computations shows that
\[
     t<(1-\frac{\epsilon}{p})\nu_{k-1}.
\]
Since $t\in\Lambda_\epsilon$ by assumption, we conclude that
$p^{k-1}t\in\Lambda$. But then
\[
     p^kt_i'=p^k+ip^{k-1}t-ip^{k-1}\epsilon \in\Lambda,
\]
which shows that $t_i'\in\Lambda_\epsilon$. If $t_i<(1-\epsilon/p)\nu_k$ then
we also have $t_i'<(1-\epsilon/p)\nu_k$, and the previous argument shows that
$p^{k-1}t\in\Lambda$. As before, we can conclude that
\[
    p^kt_i = p^k+ip^{k-1}t \in\Lambda
\]
and hence $t_i\in\Lambda_\epsilon$. Now (ii) is proved. 

Given a strictly increasing sequence $t_1<t_2<\ldots\in\Lambda_\epsilon$, we
distinguish two cases. In the first case, there exists $k_0$ such that $t_i\in
p^{-k_0}\Lambda$, for all $i$. Then the sequence is unbounded, and the claim
in (iii) is correct. In the other case, there exists a strictly increasing
sequence of indices $i_1<i_2<\ldots$ such that
$p^kt_{i_k}\not\in\Lambda$. Since we assume $t_{i_k}\in\Lambda_\epsilon$ this
means that
\[
     t_{i_k} \geq (1-\frac{\epsilon}{p})\nu_k
        \;\longrightarrow_{k\to\infty}\; \frac{p-\epsilon}{p-1}.
\]
Now (iii) follows immediately. Finally, the proof of (iv) is straightforward. 
\Endproof

\subsection{}

For the rest of this paper we shall keep the field $K$ fixed. Recall that we
have chosen certain elements $\pi_t\in k$ in the field of constants of $K$
with $v(\pi_t)=t$, for all $t\in\Lambda=v(K^\times)$, see \S \ref{pit}. It is
clear that we can extend the definition of $\pi_t$ for all $t\in \QQ$, where
$\pi_t$ lies in a fixed algebraic closure of $k$ and such that Remark
\ref{pitrem} holds.  Let $N=p^m$ be a power of $p$ and set
$\Lambda':=\frac{1}{N}\Lambda$.  Then the field
\[
      K':=K[\pi_t\mid t\in\Lambda']
\]
is a totally ramified constant extension of $K$ such that
$v((K')^\times)=\Lambda'$. Many of the following statements depend on the
choice of $K'$ (i.e.\ on $m$) and are true only for $K'/K$ (i.e.\ $m$)
sufficiently large.

We set $M':=MK'$. Since $K'/K$ is constant and totally ramified and $M/K$ is
fierce, $M'/M$ is again constant and totally ramified, with
$[M':M]=[K':K]$. 

It is easy to see that an element $a\in\OO_{M'}$ can be uniquely written in
the form
\begin{equation} \label{canform}
   a = a_1\cdot\pi_{t_1}+\ldots+a_r\cdot\pi_{t_r},
\end{equation}
with $a_i\in\OO_M^\times$, $t_i\in\Lambda'$, $0\leq t_1<\ldots<t_r$ and such
that $t_i-t_j\not\in\Lambda$, for $i\neq j$. We call \eqref{canform} the {\em
  canonical form} and the numbers $t_1,\ldots,t_r$ the {\em exponents} of $a$.

\begin{defn} \label{moderatedef}
  An element $a\in\OO_{M'}$ is called {\em moderate} if its exponents
  $t_1,\ldots,t_r$ are contained in $\Lambda_\epsilon$.
\end{defn}

It is important for us that the condition of being moderate is compatible with
enlarging the extension $K'/K$ (i.e.\ increasing $m$). 

\begin{lem} \label{moderatelem}
\begin{enumerate}
\item
  Suppose $a\in\OO_{M'}$ can be written as
  \[
     a=a_1\cdot\pi_{t_1}+\ldots+a_r\cdot\pi_{t_r},
  \]
  with $a_i\in\OO_M$ and $t_i\in\Lambda'$, $t_i\geq 0$. If
  $t_i\in\Lambda_\epsilon$ for all $i$ then $a$ is moderate. 
\item
  If $a_1,a_2\in\OO_{M'}$ is moderate, then so is $a_1\pm a_2$ and
  $a_1a_2$. In other words: the set of moderate elements is a subring of
  $\OO_{M'}$. 
\item Let $u\in\OO_{M'}^\times$ be a moderate principal unit. Then
  \[
       v(u-1)\in\Lambda_\epsilon.
  \]
  Furthermore, $u^{-1}$ is again moderate.
\item Let $u\in\OO_{M'}^\times$ be a moderate principal unit. Let $\psi_u\in
  H^1_p(M')$ denote the Kummer character associated to $u$. Assume that
  \[
     \sw(\psi_u) > \frac{\epsilon}{p-1}.
  \]
  Then (after enlarging the constant extension $K'/K$) $u\sim_{M'}u'$, where
  $u'\in\OO_{M'}^\times$ is moderate and reduced (in the sense of Definition
  \ref{reduceddef}).
\item For $t\in\Lambda_\epsilon$, with $t\geq\epsilon$, and
  $a\in\OO_K\cap\tilde{\OO}_M^{(\epsilon)}$ (see \S \ref{epsilonsubsec} for
  notation), we set
  \[
        u:=1+\pi_{t-\epsilon}\cdot a \in\OO_{M'}^\times.
  \]
  Then (after enlarging the constant extension $K'/K$) we have $u\sim_{M'}u'$,
  for a moderate principal unit $u'\in\OO_{M'}^\times$.
\end{enumerate}
\end{lem}

\proof 
Let $a=\sum_ia_i\pi_{t_i}$ be as in (i). Let $I\subset\{1,\ldots,r\}$ be the
subset of indices $i$ such that
$t_i\leq t_j$ for all $j$ with $t_i-t_j\in\Lambda$. Then 
\[
    a=\sum_{i\in I}\, a_i'\pi_{t_i},
\]
with
\[
   a_i':=\sum_{t_j-t_i\in\Lambda} a_j\pi_{t_j-t_i} \in\OO_M.
\]
We can write $a_i'=a_i''\pi_{s_i}$ with $a_i''\in\OO_M^\times$ and
$s_i\in\Lambda$. Then
\[
    a = \sum_{i\in I} \,a_i''\pi_{t_i+s_i}
\]
is, up to reordering the terms of the sum, the canonical form of $a$. Since
$t_i\in\Lambda_\epsilon$ by assumption, $t_i+s_i\in\Lambda_\epsilon$ as
well. Hence $a$ is moderate, proving (i). Assertion (ii) follows easily from (i)
and Lemma \ref{Lelem} (i).

Let $u=1+\pi_t a$ be moderate, with $a\in\OO_{M'}^\times$ and
$t>0$. We have to show that $t\in\Lambda_\epsilon$. Let
$a=a_0+a_1\pi_{t_1}+\ldots+a_r\pi_{t_r}$ be the canonical form. Then
$t_i\not\in\Lambda$, for all $i$. Applying the reasoning of the proof of (i) to
\[
    u = 1+a_0\pi_t+a_1\pi_{t_1+t}+\ldots a_r\pi_{t_r+t}
\]
we see that $t$ is either an exponent of $u$ or $t\in\Lambda$. In both
cases, $t\in\Lambda_\epsilon$, proving the first part of (iii). It remains to
show that $u^{-1}$ is also moderate. We write
\[
      u^{-1}=1-\pi_ta+\pi_{2t}a^2-\ldots = 1 +\sum_{i1}^r a_i\pi_{t_i},
\]
where the last expression is the canonical form of $u^{-1}$. Arguing again
with the proof of (i) we see that for all $i$ there exists an integer $k\geq
0$ such that $t_i=kt+s$, with $s\in\Lambda$, $s\geq 0$. Hence
$t_i\in\Lambda_\epsilon$ and $u^{-1}$ is moderate.

For the proof of (iv) we continue with the same notation. If
$\bar{a}\not\in\bar{M}^p$ then $u$ is reduced and we are done. Otherwise, we
can find $b\in\OO_M^\times$ with $\bar{b}^p=\bar{a}$ (see Definition
\ref{reduceddef}). After enlarging the extension $K'/K$ we may assume that
$t/p\in\Lambda'$. Set $v:=1+\pi_{t/p}\cdot b\in\OO_{M'}^\times$. Then
\[
    v^p=1+\pi_t\cdot b^p + \sum_{i=1}^{p-1} \pi_{1+it/p}\cdot b_i,
\]
with $b_i\in\OO_M$. Since $t\in\Lambda_\epsilon$, (i) and Remark
\ref{Lelem} (ii) show that $v^p$ is moderate. Write
\[
    u_1:=uv^{-p}=1+\pi_{t_1}\cdot a_1,
\]
with $a_1\in\OO_{M'}^\times$. Then $t_1>t$, and it follows from (i) and (iii)
that $u_1$ is moderate. If $u_1$ is reduced, we are done. Otherwise, we apply
the same procedure again to $u_1$. Continuing this way, we obtain a sequence
of moderate principal units $u_1,u_2,\ldots$ with $u_i\sim_{M'} u$ and
$u_i=1+\pi_{t_i}\cdot a_i$ such that
\begin{equation} \label{modlemeq2}
   t<t_1<t_2<\ldots .
\end{equation}
We have to show that this process stops after a finite number of steps with a
reduced principal unit $u_k$, for some $k\in\NN$. We argue by contradiction
and assume that we obtain an infinite sequence of $u_i$. 
Since $u_i$ is moderate, $t_i\in\Lambda_\epsilon$ for all $i$, by
(iii). Therefore, Lemma \ref{Lelem} (iii) shows that
\begin{equation} \label{modlemeq3}
    \limsup_i t_i \geq \frac{p-\epsilon}{p-1}.
\end{equation}
On the other hand we have
\begin{equation} \label{modlemeq4}
  t_i\leq \frac{p}{p-1}-\sw(\psi_u) < \frac{p-\epsilon}{p-1}
\end{equation}
by Proposition \ref{reducedprop2} and the assumption. But \eqref{modlemeq3}
and \eqref{modlemeq4} contradict each other. This completes the proof of (iv).

It remains to prove (v). By assumption we can write
$a=b^p+\pi_\epsilon\cdot c$, with $b,c\in\OO_M$, and hence
\[
     u = 1 +\pi_{t-\epsilon}\cdot b^p+\pi_t\cdot c = u_0+u_1,
\]
with $u_0:=1 +\pi_{t-\epsilon}\cdot b^p$ and $u_1:=\pi_t\cdot c$.
After enlarging $K'/K$ we may assume that $(t-\epsilon)/p\in\Lambda'$. We set
$w:=1+\pi_{(t-\epsilon)/p}\cdot b\in\OO_{M'}$. Then
\[
    w^p=1+\pi_{t-\epsilon}\cdot b^p + 
         \sum_{i=1}^{p-1}\pi_{1+i(t-\epsilon)/p}\cdot b_i
       = u_0 + w_1,
\]
with $b_i\in\OO_M$ and $w_1\in\OO_{M'}$. We have $v(w_1)>1>\epsilon$ and,
moreover, $w_1$ is moderate by Lemma \ref{Lelem} (ii). Using (ii) and Lemma
\ref{Lelem} (v) one shows that
\[
   u_0^{-1}w_1 =(1-\pi_{t-\epsilon}b^p+\pi_{2(t-\epsilon)}b^{2p}-\ldots)w_1
\]
and hence 
\[
      w_2:=(1+u_0^{-1}w_1)^{-1}
\]
is moderate. A similar argument shows that $u_0^{-1}u_1$ is moderate.  
Using once more (ii) we conclude that
\[
   u':=uw^{-p}=(1+u_1u_0^{-1})w_2
\]
is moderate. This finishes the proof of the lemma.
\Endproof

\subsection{}

We can now give a proof of Lemma \ref{minlem}. We argue by
contradiction, i.e.\ we assume that $\sw(\chi)>p\deltab$ for all lifts $\chi$
of $\chib$. Using this assumption, we shall construct inductively a sequence
of lifts $\chi_0,\chi_1,\chi_2,\ldots$ of $\chib$ whose Swan conductors
$\delta_i:=\sw(\chi_i)$ form a strictly decreasing sequence. Afterwards we
will show that $\lim_i\delta_i< p\deltab$, which gives the desired
contradiction. For a fixed $i$, the lifts $\chi_0,\ldots,\chi_i$ will be
defined and fierce over a finite constant extension $K'/K$ as described
above. During the induction process, we will have to keep increasing $K'/K$,
but this is not a problem.

We start by choosing an arbitrary lift $\chi_0\in H^1_{p^n}(K)$ of $\chib$
defined over $K$. We may assume that $\chi_0$ is fierce over $K'$. By
assumption we have $\delta_0:=\sw(\chi_0)>p\deltab$. By induction, suppose
that $\chi_i$ has been constructed and is defined and fierce over the constant
extension $K'/K$. Set $\delta_i:=\sw(\chi_i)$ and $\omega_i:=\dsw(\chi_i)$. We
have $\deltab<1/(p-1)$ and $p\deltab<\delta_i$ by assumption. Therefore,
Proposition \ref{indprop1} says that $\delta_i\leq p/(p-1)$ and that
$\C(\omega_i)=0$ if $\delta<p/(p-1)$ and $\C(\omega_i)=\omega_i$
otherwise. Note that the second case, $\delta_i=p/(p-1)$, can occur only once,
for $i=0$. It follows that $\omega_i=d\bar{a}_i$ in the first and
$\omega_i=d\bar{a}_i/\bar{a}_i$ in the second case, for an element
$\bar{a}_i\in\bar{K}$. Using Corollary \ref{epsiloncor} we can choose
$a_i\in\OO_K\cap\tilde{\OO}_M^{(\epsilon)}$ lifting $\bar{a}_i$. Set
$\mu_i:=p/(p-1)-\delta_i$ and
\[
    z_i:=1- \pi_{\mu_i}\cdot a_i\in\OO_{K'}^\times
\]
if $\delta_i<p/(p-1)$ (and $z_i:=a_i^{-1}$ otherwise). Here we assume that
$\mu_i\in\Lambda'$. Then, as an element of $\OO_{K'}$, $z_i$ is reduced by
Proposition \ref{reducedprop}. Let $\psi_i\in H^1_p(K')$ be the character
associated to $z_i$ by Kummer theory. By Proposition \ref{reducedprop2}
$\psi_i$ is fierce over $K'$ and we have
\[
   \sw(\psi_i)=p/(p-1)-\mu_i=\delta_i, \quad 
       \dsw(\psi_i)=-\omega_i.
\]
We define $\chi_{i+1}:=\chi_i+\psi_i\in H^1_{p^n}(K')$. This is again a lift
of $\chib$, and Lemma \ref{rswlem} shows that
$\delta_{i+1}:=\sw(\chi_{i+1})<\delta_i$. This completes the construction of
the sequence $\chi_0,\chi_1,\ldots$. 

The character $\chi_0$ corresponds, via Kummer theory, to the class of an
element $u\in K^\times$ in $K^\times/(K^\times)^{p^n}$. It follows that
$M=K[y]$, where $y^{p^{n-1}}=u$. Moreover, the restriction $\chi_0|M\in
H^1_p(M)$ corresponds to the class of $y$ in $M^\times/(M^\times)^p$.

Let $\chi\in H^1_{p^n}(K')$ be an arbitrary lift of $\chib$ defined over
$K'$. Then $\chi|_{M'}$ corresponds to the class of $yz$ in
${M'}^\times/({M'}^\times)^p$, for some $z\in {K'}^\times$. Note that $\chi$
is fierce if and only if $\chi|_{M'}$ is, and this is the case if and only if
$yz$ is $M'$-equivalent to a principal unit of the form $y'=1+\pi_\lambda\cdot
w$, with $w\in\OO_{M'}^\times$ and $\bar{w}\not\in\bar{M}^p$. Assume that this
is the case. Then $\deltat:=\sw(\chi|_{M'})=p/(p-1)-\lambda$. Using
\eqref{psieq} we obtain
\begin{equation} \label{deltalambdaeq}
   \delta= \frac{p}{p-1} +\epsilon -\lambda.
\end{equation}

For the sequence of lifts $\chi_0,\chi_1,\ldots$ constructed above, this means
the following. The restriction $\chi_i|_{M'}\in H^1_p(M')$ corresponds to the
class of $y_i\in{M'}^\times$, defined by the recursive formula $y_0:=y$,
$y_{i+1}:=y_iz_i$. 

\begin{lem} \label{claimlem}
 We have $y_i\sim_{M'}y_i'$, where $y_i'\in\OO_{M'}^\times$ is 
 reduced and moderate.
\end{lem}

\proof
We prove the lemma by induction on $i$, using Lemma \ref{moderatelem}
repeatedly. For $i=0$ there is nothing to show. By induction, we may assume
that $y_i\sim_{M'}y_i'$, with $y_i'$ reduced and moderate. Write
\[
   y_i'=1+\pi_{\lambda_i}\cdot v_i,
\]
with $\bar{v}_i\not\in\bar{M}^p$. Then \eqref{deltalambdaeq} says that
$\delta_i=\sw(\chi_i)=p/(p-1)+\epsilon-\lambda_i$. Therefore,
\[
     \mu_i=p/(p-1)-\delta_i=\lambda_i-\epsilon.
\]
Since $\lambda_i\in\Lambda_\epsilon$ by Assertion (iii) of Lemma
\ref{moderatelem}, Assertion (v) of that lemma shows that $z_i$ is moderate
(here we use that $a_i\in\tilde{\OO}_M^{(\epsilon)}$). Hence $y_i'z_i$ is
moderate by Assertion (ii).  The Kummer class of $y_i'z_i$ is the character
$\chi_{i+1}|_{M'}$. We have
\[\begin{split}
   \sw(\chi_{i+1}|_{M'})-\frac{\epsilon}{p-1}
          & =\delta_{i+1}-\frac{p}{p-1}\cdot\epsilon\\
    & >p\deltab-\big(\delta_{n-1}+\frac{1}{p}\delta_{n-2}+\ldots
                +\frac{1}{p^{n-1}}\delta_1\big) \\
    & \geq p\deltab-\deltab\big(1+\frac{1}{p^2}+\ldots
                +\frac{1}{p^{2n-2}}\big) \\
    & >\deltab\cdot\big(p-\frac{p^2}{p^2-1}\big) > 0.
\end{split}\]
So Assertion (iv) of Lemma \ref{moderatelem} shows that
$y_{i+1}\sim_{M'}y_i'z_i\sim_{M'}y_{i+1}$, with $y_{i+1}'$ reduced and
moderate. This completes the proof of Lemma \ref{claimlem}. 
\Endproof

Since $y_i'$ is reduced by Lemma \ref{claimlem}, \eqref{deltalambdaeq} shows
that
\[
    \lambda_i=v(y_i'-1) = p/(p-1)+\epsilon-\delta_i
\]
It follows that the sequence $\lambda_0<\lambda_1<\ldots$ is strictly
increasing. On the other hand, $y_i'$ is also moderate, so
$\lambda_i\in\Lambda_\epsilon$ by Lemma \ref{moderatelem} (iii). Therefore,
Lemma \ref{Lelem} (iii) implies
\[
    \limsup_i \lambda_i \geq \frac{p-\epsilon}{p-1}.
\]
We conclude that 
\begin{equation} \label{mineq10}
   \liminf_i \delta_i \leq \frac{p}{p-1}+\epsilon-\frac{p-\epsilon}{p-1}
        = \frac{p}{p-1}\cdot\epsilon.
\end{equation}
By the calculation already used in the proof of Claim 3 we have
\begin{equation} \label{mineq11}
   \frac{p}{p-1}\cdot\epsilon < \deltab\cdot\frac{p^2}{p^2-1} < p\,\deltab.
\end{equation}
Combining \eqref{mineq10} and \eqref{mineq11} shows that $\delta_i<p\deltab$
for $i$ sufficiently large. But this contradicts Proposition \ref{indprop1},
and Lemma \ref{minlem} follows. 
\Endproof

\end{document}